\documentclass[final]{siamltex}

\usepackage{epsf,latexsym,amsfonts,amssymb,subfigure,mathrsfs,graphicx,stmaryrd,
amsmath}
\usepackage[notref,notcite]{showkeys}
\usepackage{color}

\numberwithin{equation}{section}

\newcommand{\Lr}[1]{\left(#1\right)}
\newcommand{\mc}[1]{\mathcal #1}
\newcommand{\mb}[1]{\mathbb #1}

\newcommand{\ov}[1]{\overline{#1}}
\newcommand{\abs}[1]{\lvert#1\rvert}

\def\al{\alpha}
\def\ga{\gamma}
\def\ka{\kappa}
\def\del{\delta}
\newcommand{\nn}{\nonumber}

\title{A study on the quasiconinuum approximations of a one-dimensional fracture model}
\author{Xiantao Li
\thanks{Department of Mathematics, the Pennsylvania State University, University Park,
Pennsylvania, 16802, U.S.A. ({\tt xli@math.psu.edu})
{The work of Li was supported by National Natural Science
              Foundation grant DMS1016582.}}
\and{Pingbing Ming}
\thanks{LSEC, Institute of Computational
  Mathematics and Scientific/Engineering Computing,
  AMSS, Chinese Academy of Sciences,
  No. 55, Zhong-Guan-Cun East Road,
  Beijing 100190, China. ({\tt mpb@lsec.cc.ac.cn}) {The work of Ming was supported by National Natural Science Foundation of China grants 10932011 and 91230203, and by the funds from Creative Research Groups of China through grant 11021101, and by the support of CAS National Center for Mathematics and Interdisciplinary Sciences.}}}
\date{\today}

\begin{document}
\maketitle
\begin{abstract}
We study three quasicontinuum approximations of a lattice model for crack
propagation. The influence of the approximation on the bifurcation patterns is investigated.  The estimate of the modeling error is applicable to near and beyond bifurcation points, which enables us to evaluate the approximation over a finite range of loading and multiple mechanical equilibria.
\end{abstract}


\section{Introduction}
In recent years multiscale models have undoubtedly become one of the most important computational tools for problems in materials science.  Such multiscale models allow atomistic details of local defects, while taking advantage of the efficiency of continuum models  to handle the  calculations in the majority of the computational domain. One remarkable success in multiscale modeling of materials science is the quasi-continuum (QC) method \cite{TadmorOrtizPhillips:1996}, which couples a molecular mechanics
model with a continuum finite element model. The QC method has motivated a lot of recent works on multiscale models of crystalline solids \cite{KnapOrtiz:2001,BeXi03,WaLi03,Shimokawa:2004,LI12a}.

Meanwhile, there has been considerable interest from the applied mathematics community to analyze the stability and accuracy of QC type methods \cite{Lin03,EMing:2005, ELuYang:2006,Ming:2008,MingYang:2009,DobsonLuskin:2009a,DobsonLuskin:2009b,DobsonLuskinOrtner:2010a,Shapeev:2011,ChenMing:2012, OrtnerShapeev,Ortner,CuiMing:2013, LuMing:2013,LiMing:2013}. Various important issues, such as the ghost forces and stability, have been extensively investigated.  One major weakness of all the existing results, however, is that they are only applicable to a system near {\it one} local minimum with a fixed load, even
for a relatively simple system~\cite{Jusuf2010}. This significantly limits the practical values of these analysis. First, for any given loading condition, there are typically a large number of local mechanical equilibria. Second, often of interest in practice is the transition of the system as the loading condition changes. Examples include phase transformation \cite{ElShTr02}, crack propagation and kinking \cite{ChYi94,CoRi80} and dislocation nucleation \cite{PlCaBo08}, etc. Throughout these processes, the system is driven from a stable equilibrium to a critical point, where the system loses its stability, and then settles to another equilibrium.

Fortunately, the theory of bifurcation \cite{Sotomayor73,Carr81,GuHo83,Benoit90,ChowHale,MaWa05} provides a rigorous tool to understand the transition processes. The theory considers models, either static or dynamic,  with certain embedded parameters, which for mechanics problems, naturally correspond to the external loading conditions. Bifurcation arises when the system loses its stability,
and it is a ubiquitous phenomenon in mechanics \cite{MarHug94}. A reduction procedure is available \cite{Carr81} to probe the transition process. Of particular interest in this context is the bifurcation diagram, consisting of bifurcation curves for a wide range of parameters. The curves contains local equilibria, including both stable and unstable ones. As a result, the analysis is well beyond local, stable equilibrium. This is the primary motivation for the current work.

The molecular mechanics model becomes highly indefinite at the bifurcation point, and
the standard analysis is not applicable due to the loss of coercivity.
 In fact, most existing results relied on even more strict stability conditions. We refer  to~\cite{EMing:2005, MingYang:2009,DobsonLuskin:2009b,DobsonLuskinOrtner:2010a} for related discussion. There are some methods that have sharp stability conditions~\cite{LuMing:2013}.
Nevertheless, these analysis do not predict the modeling error beyond the bifurcation point.

As a first attempt to understand the modeling error of atomistic-to-continuum coupling methods over a more global scale, we consider the QC approximations near
and beyond bifurcation points, when applied to a one-dimensional fracture model. The first of such models has been constructed  in   \cite{ThHsRa71,FuTh78} to understand the atomic aspect of crack initiation, which led to the important concept of lattice trapping. We have modified the original model so that the QC methods can be directly
applied. In particular, three methods, including the original QC method, the quasi-nonlocal QC method \cite{Shimokawa:2004},  and a force-based method~\cite{KnapOrtiz:2001}, are considered in this paper. For each method, we derive an effective equation that describes the bifurcation diagram.
This is in the same spirit as the centre manifold~\cite{Carr81}, a tool that significantly reduces the dimension of the problem.
The one-dimensional lattice model, despite its simplicity, gives rise to bifurcation patterns that resemble those of high dimensional fracture models~\cite{XLi2013A}. Therefore, it already captures the essential mechanism behind crack initiation.

This provides a new approach to measure the modeling error: Instead of comparing the atomic displacement, which may not have an error bound near bifurcation points, we compare the bifurcation curves. Intuitively, when the bifurcation curves are accurately produced, the transition mechanism is well captured.
To quantitatively estimate the error in predicting the bifurcation curves, we formulate the bifurcation equations  as solutions of some ordinary differential equations. Then, the difference between the bifurcation curves for the full atomistic model and the coupled models can be estimated using stability theory of ordinary differential equations. Since this is a new issue that has not been addressed in previous works,
we have chosen the simple lattice model of fracture to illustrate the ideas. For this particular example, we are able to find explicitly the parameters in the bifurcation diagram, and make direct comparisons. The extension to more general problems will be investigated in future works.

The rest of the paper has been organized as follows. In \S~\ref{sec: crack}, we introduce the lattice model and find the explicit solution of this model.
The bifurcation behavior is also discussed. In \S~\ref{sec: multi}, we obtain the bifurcation diagrams of three QC approximations.
In \S~\ref{sec: ode}, we quantify the difference among the bifurcation curves.
\section{The Lattice Fracture Model}\label{sec: crack}
\subsection{The lattice model}

\begin{figure}[htbp]
\begin{center}
\includegraphics[scale=0.5]{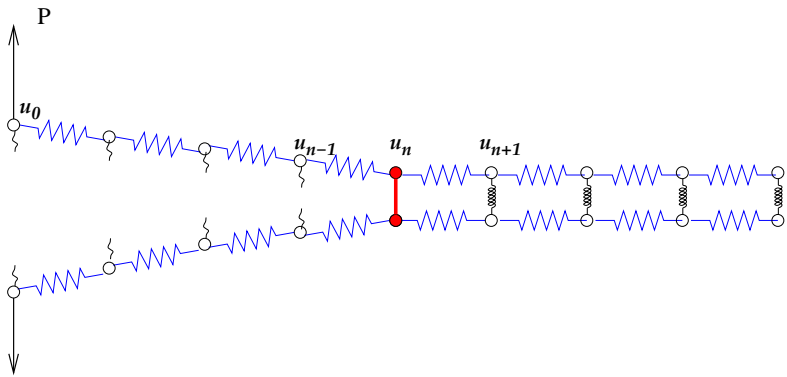}
\caption{A schematic of the lattice model. The springs indicate linear interactions
between two atoms. The solid line represents a nonlinear bond at the crack-tip.  }
\label{fig: lattice}
\end{center}
\end{figure}

We consider a one-dimensional chain model for crack propagation, which has been used to study the lattice trapping effect~\cite{ThHsRa71,FuTh78}. The system consists of two chains of atoms above and below the crack face.  The atoms are only allowed to move vertically. The displacement of the atoms in the chain below the crack face are assumed to be opposite to the corresponding atoms in the upper chain. Hence only the atoms in the upper chain need to be considered. This is illustrated in Fig.~\ref{fig: lattice}.

Each atom of the chain is interacted with the two nearest atoms on the left and the two nearest atoms on the right.  In addition, it is interacted with the atom below (or above) via a nonlinear force, which is denoted by $F(u)$ and satisfies the following conditions:
\begin{enumerate}
\item $\ka_3=-F'(0)>0$;
\item $F(u)=0,$ if $u>u_{\rm cut}$;
\item $F'(u_{\rm cut})=0$.
\end{enumerate}
Here \(u_{\rm cut}\) is a cut-off distance, and bonds are considered to be {\it broken} beyond this threshold. The last condition is a smoothness assumption often made in the  analysis. The following simple example of $F(u)$ satisfies all those conditions,
\begin{equation}\label{eq: F(u)}
  F(u) = -\frac{\ka_3}{u_{\rm cut}^2}u(u- u_{\rm cut})^2 \chi_{[0, u_{\rm cut}]}(u),
\end{equation}
where $\chi_{[0, u_{\rm cut}]}$ is an indicator function.

We assume that the vertical bonds are already broken for \(j<n\) with $n$
the crack-tip position. This creates an existing crack and allows us to study crack propagation.
We further simplify the model by replacing the nonlinear bonds ahead of the crack-tip by
linear springs with spring constant $\ka_3$.

The surface energy density is defined by
 \[
 \gamma(u){:}= -\int_0^u F(v)dv,
 \]
 and we denote \(\gamma_0{:}= \gamma(u_{\rm cut}).\)

The system is subject to a loading \(P\) at the left most atom.
The total energy for the upper chain reads as
\[
 \begin{aligned}
  E&=-Pu_0 + \sum_{j\ge 0}\Lr{
  \dfrac{\ka_1}{2}\big(u_{j+1}-u_j\big)^2
  +\dfrac{\ka_2}{2}\big(u_{j+2}-u_j\big)^2}\\
  &\quad+ 2n \gamma_0 + 2\gamma(u_n) + {\ka_3} \sum_{j>n} u_j^2,
  \end{aligned}
\]
where \(\ka_1\) and \(\ka_2\) are respectively the force constants for the
nearest and next nearest neighbor interactions. We assume the force constants satisfy
\begin{equation}\label{eq: k-cond}
 \ka_1>0, \; \bar{\ka}{:}=\ka_1+4\ka_2>0, \; \ka_3=\gamma''(0)>0.
\end{equation}

We denote the cracked region as \(A\), and define
\[
L^A u_j {:}=\ka_1(u_{j+1}-2u_j+u_{j-1}) + \ka_2(u_{j+2}-2u_j+u_{j-2}) .
\]
Similarly, we denote the un-cracked region be \(B\) and define
\(
L^B u_j {:}=L^A u_j-2\ka_3 u_j\).
The force balance equations are given by
\begin{align}
\ka_1(u_1-u_0) + \ka_2(u_2-u_0) + P&=0, \label{eq:leftbc}\\
\ka_1(u_2-2u_1+u_0) + \ka_2(u_3-u_1) &=0,\label{eq:leftbc1}\\
L^A u_j&=0,\quad j=2,\dots,n-1,\label{eq: cracked}\\
L^A u_n+F(u_n)&=0,\label{eq: tip}\\
L^B u_j&=0,\quad j\ge n+1.\label{eq: un-cracked}
\end{align}
\subsection{The solution near the crack tip}
In this section, we study the solution at the crack tip by eliminating other degrees of freedom.
We start with the atoms along the crack face, where we have a difference equation
with the following characteristic equation
\begin{equation}\label{eq: pA}
p^A(z)=0, \quad p^A(z){:}= \ka_2z^4 + \ka_1 z^3 -(2\ka_1+2\ka_2)z^2 + \ka_1z + \ka_2.
\end{equation}
We factor $p^A(z)$ as
\(
p^A(z)=(z-1)^2(z-z_0)(z-z_0^{-1}),
\)
where
\[
z_0=-1-\dfrac{\ka_1}{2\ka_2}(1-\sqrt{\ov{\ka}/\ka_1}\,).
\]

By~\eqref{eq: k-cond}, one can verify that $z_0<1$, and \(z_0\) and \(1/z_0\) solve
\begin{equation}\label{eq: z0}
  \ka_2 z^2 + (\ka_1 + 2\ka_2)z + \ka_2=0.
\end{equation}

Next we turn to the region ahead of the crack tip, where the characteristic equation is
\begin{equation}\label{eq: pB}
p^B(z)=0, \quad p^B(z){:}= \ka_2z^4 + \ka_1 z^3 -(2\ka_1+2\ka_2+2\ka_3)z^2 + \ka_1z + \ka_2.
\end{equation}
In this case, the general solutions can be written as
\begin{equation}\label{eq: uB}
  u^B_j= B_1 z_1^j + B_2 z_2^j,
\end{equation}
where \(z_1\) and \(z_2\) are two roots of the characteristic equation that are less or equal to one. We focus on the case when all the roots are real. This occurs when
\(
\Delta= \ov{\ka}^2+8\ka_2\ka_3\ge 0.
\)

Once we have \(z_1\) and \(z_2\), the polynomial can be factored into
\[
p^B(z) = \ka_2 (z-z_1)(z-z_1^{-1})(z-z_2)(z-z_2^{-1}).
\]
By comparing the coefficients, we find 
\begin{equation}\label{eq: z3}
\ka_1 z_1 z_2=-\ka_2 (1+z_1z_2)(z_1 + z_2).
\end{equation}
This equation will be used later to simplify our calculation.

At the interface, we have the matching conditions:
\begin{equation}\label{eq:match}
u^A_i=u^B_i\qquad i=n,\;n-1.
\end{equation}
For brevity, we drop the superscripts \(A\) and \(B\) whenever there is no confusion. By setting \(j\) to
\(n-1, n, n+1\) and \(n+2\) in equation \eqref{eq: uB}, we find
\[
u_{n+1} =\alpha u_{n-1} + \beta u_n,
\]
and
\[
u_{n+2} =\al u_n+\beta u_{n+1}= \alpha \beta u_{n-1} + (\alpha+\beta^2) u_n,
\]
where
\(
\alpha=- z_1 z_2,\;\beta=z_1 + z_2.
\)
 For other atoms in this region, the displacement can be obtained recursively as
 \[
 u_{j+1}= \alpha u_{j-1} + \beta u_j\quad\text{for any\;}j \ge n+1.
 \]
 In terms of the strain, these conditions can be expressed as
\begin{equation}\label{eq: matching-tip}
\begin{aligned}
  u_{n+1} - u_n & = - \alpha (u_n - u_{n-1}) + (\al+\beta - 1) u_n,\\
  u_{n+2} - u_n & = - \alpha\beta (u_n - u_{n-1}) + (\alpha\beta +\beta^2 +\alpha - 1) u_n.
\end{aligned}
\end{equation}

By~\eqref{eq: z3}, we find
\begin{equation}\label{eq:alphabeta}
  \ka_1 \alpha + \ka_2(\alpha-1) \beta=0, \quad 2(\ka_1+\ka_2)
  \alpha  + \ka_2(\alpha^2 + \beta^2 + 1)= 2\ka_3 \alpha.
\end{equation}

With these preparation, we are ready to find solutions in the crack region. For $j\le n+1$, we express the solution as
\begin{equation}\label{eq:uA}
u_j=a+bj+c\cosh[j\del]+d\sinh[j\del]
\end{equation}
with
\begin{equation}\label{eq:pararela1}
\cosh\del=-1-\ka_1/(2\ka_2).
\end{equation}
In particular, we choose \(\del= -\log z_0.\)

We proceed to derive an equation for $u_n$ by eliminating
all other variables in $L^A u_n$. It follows from ~\eqref{eq: matching-tip} that
 \begin{align}
L^A u_n
  &=  \ka_2(u_{n+2}- u_n) + \ka_1(u_{n+1}-u_n)\\
&\quad- (\ka_1+\ka_2)(u_n - u_{n-1}) - \ka_2(u_{n-1}-u_{n-2})\nn\\
 &=
   [\ka_1(\alpha+\beta-1)+\ka_2(\alpha\beta+\beta^2+\alpha-1)\big] u_n\nn\\
   &\quad-\big[\ka_1(1+\alpha)+\ka_2(1+\alpha\beta)\big](u_n-u_{n-1})
   -\ka_2(u_{n-1} - u_{n-2})\nn\\
&=(\alpha+\beta-1)\bigl(\ka_1+\ka_2(\beta+1)\bigr)u_n\nn\\
   &-\Lr{\bigl(\ka_1+\ka_2(1+\beta)\bigr)(u_n-u_{n-1})+\ka_2(u_{n-1} - u_{n-2})},\label{eq:mach-exp}
\end{align}
where we have used the identity
\[
\ka_1(1+\alpha)+\ka_2(1+\alpha\beta)=\ka_1+\ka_2(1+\beta),
\]
which follows from~\eqref{eq:alphabeta}. To calculate 
the second term in~\eqref{eq:mach-exp},
we use the following relations that can be easily verified, and the proof
can be found in the Appendix. 

For any $k\in\mb{Z}$, there holds
\begin{align}
(\ka_1+2\ka_2)(\cosh[k\del]-\cosh[(k-1)\del])&+\ka_2(\cosh[(k-1)\del]-\cosh[(k-2)\del])\nn\\
&=-\ka_2(\cosh[(k+1)\del]-\cosh[k\del]),\label{eq:rela1}
\end{align}
and
\begin{align}
(\ka_1+2\ka_2)(\sinh[k\del]-\sinh[(k-1)\del])&+\ka_2(\sinh[(k-1)\del]-\sinh[(k-2)\del])\nn\\
&=-\ka_2(\sinh[(k+1)\del]-\sinh[k\del]).\label{eq:rela2}
\end{align}

For any $k\in\mb{Z}$ and $\rho\in\mb{R}$, we define
\begin{align*}
\mc{F}_{k,\rho}(\del)&{:}=\cosh[(k+1)\del]-(1-\rho)\cosh[k\del]-\rho\cosh[(k-1)\del],\\
\mc{G}_{k,\rho}(\del)&{:}=\sinh[(k+1)\del]-(1-\rho)\sinh[k\del]-\rho\sinh[(k-1)\del].
\end{align*}
Using~\eqref{eq:rela1} and~\eqref{eq:rela2} with $k=n$ and $\rho=1-\beta$, we obtain
\begin{align*}
\bigl(\ka_1+\ka_2(\beta+1)\bigr)(u_n-u_{n-1})
&+\ka_2(u_{n-1}-u_{n-2})\\
&=\bigl(\ov{\ka}+\ka_2(\beta-2)\bigr)b
-\ka_2\bigl(\mc{F}_{n,1-\beta}(\del)c+\mc{G}_{n,1-\beta}(\del)d\bigr).
\end{align*}
Substituting the above identity into~\eqref{eq:mach-exp}, we obtain
\begin{equation}\label{eq:match-expa}
 \begin{aligned}
L^A u_n&=(\alpha+\beta-1)\bigl(\ka_1+\ka_2(\beta+1)\bigr)u_n\\
   &\quad-\bigl(\ov{\ka}+\ka_2(\beta-2)\bigr)b
   +\ka_2\bigl(\mc{F}_{n,1-\beta}(\del)c+\mc{G}_{n,1-\beta}(\del)d\bigr).
\end{aligned}
\end{equation}

It remains to find the parameters $b$, $c$ and $d$. First we substitute the expressions for $u_{n+1}-u_n$ and $u_n-u_{n-1}$ into~\eqref{eq: matching-tip} and obtain
\begin{equation}\label{eq:matchcoef}
\mc{F}_{n,\al}(\del)c+\mc{G}_{n,\al}(\del)d=-(1+\al)b+(\al+\beta-1)u_n.
\end{equation}

Next we shall use the equations for $j=1,2$ to determine two parameters in
$u_j$. A simple trick is to introduce one more atom to the left, with displacement,
$u_{\bar{1}}$, and extends the equation to $j=1$,
\[
\ka_1(u_2-2u_1+u_0)+\ka_2(u_3-2u_1+u_{\bar1})=0,
\]
which together with~\eqref{eq:leftbc1} leads to \(u_1=u_{\bar1}\).
This immediately implies
\[
b=-d\sinh\del.
\]
Substituting the expression of $u_j$ into~\eqref{eq:leftbc}, we obtain
\[
(\cosh\del-1)\Lr{\ka_1+2\ka_2(\cosh\del+1)}c+2\ka_2\sinh\del(\cosh\del-1)d+P=0.
\]
Using~\eqref{eq:pararela1}, we obtain
\[
d=-\dfrac{P}{2\ka_2\sinh\del(\cosh\del-1)}=\dfrac{P/\ov{\ka}}{\sinh\del},
\]
which in turn implies \(b=-P/\ov{\ka}.\)
Using~\eqref{eq:matchcoef},\footnote{This relation is possible because we
have assumed that only the roots with magnitude less than one in the expression
of $u_j$.}, we obtain
\[
c=\dfrac{\al+\beta-1}{\mc{F}_{n,\al}(\del)}u_n
-\dfrac{P/\bar\ka}{\mc{F}_{n,\al}(\del)}\Lr{\dfrac{\mc{G}_{n,\al}(\del)}{\sinh\del}
-(1+\al)}.
\]
Substituting the expressions of $b,c$ and $d$ into~\eqref{eq:match-expa}, we obtain
an equation for $u_n$:
\begin{equation}\label{eq: eff0}
F(u_n)+\ka u_n+\eta P=0
\end{equation}
with
\[
\ka=(\alpha+\beta-1)\Lr{\ka_1+\ka_2(\beta+1)+\ka_2\dfrac{\mc{F}_{n,1-\beta}(\del)}
{\mc{F}_{n,\al}(\del)}},
\]
and
\[
\eta=\dfrac{\ka_2}{\ov{\ka}}\Lr{\dfrac{\mc{G}_{n,1-\beta}\mc{F}_{n,\al}
-\mc{F}_{n,1-\beta}\mc{G}_{n,\al}}{\sinh\del\mc{F}_{n,\al}}
+(1+\al)\dfrac{\mc{F}_{n,1-\beta}}{\mc{F}_{n,\al}}}
+\dfrac{\ka_1+\ka_2(1+\beta)}{\ov{\ka}}.
\]

The equation~\eqref{eq: eff0} is called {\em
the effective equation} because all other degrees of freedom have been removed.
Of particular interest is the limits of $\ka$ and $\eta$ when \(n\) is large. To this ends, we write
\[
\bigl(\mc{F}_{n,\al}(\del),\mc{G}_{n,\al}(\del)\bigr)=\bigl(A_\al(\del),B_\al(\del)\bigr)
\mc{K}_n
\]
with
\(
A_\al(\del)=(1-\al)(\cosh\del-1)\) and \(B_\al(\del)=(1+\al)\sinh\del\),
and the $2$ by $2$ matrix $\mc{K}_n$ is defined by
\[
\mc{K}_n{:}=\begin{pmatrix}
\cosh[n\del]&\sinh[n\del]\\
\sinh[n\del]&\cosh[n\del]
\end{pmatrix}.
\]
A direct calculation gives
\(\ka\to\ka_0\) as $n\to\infty$ with
\[
\ka_0=(\alpha+\beta-1)\Lr{\ka_1+\ka_2(1+\beta)+\ka_2
\dfrac{A_{1-\beta}(\del)+B_{1-\beta}(\del)}
{A_\al(\del)+B_{\al}(\del)}}.
\]
This is the limit when the length of the crack reaches a macroscopic size.
In particular, we have an expansion of
$\ka$ as
\[
\ka=\ka_0-\dfrac{(\al+\beta-1)\ov{\ka}\sinh\del}{(A_\al+B_\al)^2}z_0^{2n}
+\mc{O}(z_0^{4n}).
\]

To calculate the limit of $\eta$, we write
\begin{align}
\mc{G}_{n,1-\beta}\mc{F}_{n,\al}
-\mc{F}_{n,1-\beta}\mc{G}_{n,\al}&=\det\begin{pmatrix}
\mc{F}_{n,\al}&\mc{F}_{n,1-\beta}\\
\mc{G}_{n,\al}&\mc{G}_{n,1-\beta}
\end{pmatrix}
=\det\mc{K}_n\det\begin{pmatrix}
A_\al&A_{1-\beta}\\
B_{\al}&B_{1-\beta}
\end{pmatrix}\nn\\
&=-2(\al+\beta-1)(\cosh\del-1)\sinh\del.\label{eq:crisscross}
\end{align}
The right-hand side is independent of $n$, which will be exploited later on.

Substituting the above identity into the expression of $\eta$, we obtain
\[
\eta=1+\dfrac{\ka_2}{\ov{\ka}}
\dfrac{(\beta-2)\mc{F}_{n,\al}+(1+\al)\mc{F}_{n,1-\beta}}{\mc{F}_{n,\al}}
-\dfrac{2\ka_2}{\ov{\ka}}\dfrac{(\al+\beta-1)(\cosh\del-1)}{\mc{F}_{n,\al}}.
\]
A direct calculation gives
\begin{align}
(\beta-2)\mc{F}_{n,\al}+(1+\al)\mc{F}_{n,1-\beta}
&=\Lr{(\beta-2)A_\al+(1+\al)A_{1-\beta}}
\cosh[n\del]\nn\\
&=2(\al+\beta-1)(\cosh\del-1)\cosh[n\del].\label{eq:idencoeff}
\end{align}
Using the above identity, we rewrite $\eta$ as
\begin{align}
\eta
&=1+\dfrac{2\ka_2}{\ov{\ka}}
\dfrac{(\al+\beta-1)(\cosh\del-1)(\cosh[n\del]-1)}{\mc{F}_{n,\al}}\nn\\
&=1-\dfrac{(\al+\beta-1)(\cosh[n\del]-1)}{\mc{F}_{n,\al}}.\label{eq:eta-exactfinal}
\end{align}
Letting $n$ go to infinity, we obtain $\eta\to\eta_0$ with
\[
\eta_0=1-\dfrac{\al+\beta-1}{A_\al+B_\al}.
\]
We also have the following expansion for $\eta$:
\[
\eta=\eta_0+\dfrac{2(\al+\beta-1)}{A_\al+B_\al}z_0^n+\mc{O}(z_0^{2n}).
\]

Notice that we have $\ka_0 <0$, since
\begin{equation}\label{eq:rootrela}
\al+\beta-1=-(1-z_1)(1-z_2)<0,
\end{equation}
and similarly, we have $\eta_0>1$.
\subsection{Bifurcation behavior}\label{sec: bif}
\begin{figure}[htbp]
\begin{center}
\includegraphics[scale=0.45]{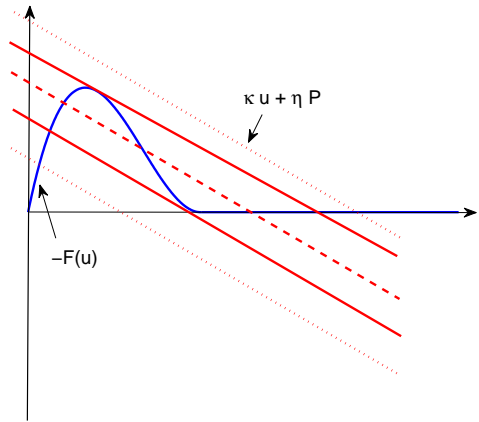}
\caption{The solutions of equation \eqref{eq: eff0}, shown as the intersections
of the function $-F(u)$ and the linear function $\kappa u + \eta P$. Dotted lines: Only one solution exists;  The dashed line: There are three solutions; Solid lines: Two of the three solutions reduce to a repeated root.   }
\label{fig: interscept}
\end{center}
\end{figure}

To understand the roles of the parameters $\kappa$ and $\eta$, we rewrite
the reduced equation~\eqref{eq: eff0} as
\begin{equation}\label{eq: eff1}
 \kappa u + \eta P= -F(u).
\end{equation}
We shall regard $\kappa$ as an intrinsic material parameter, and $P$
as an external load that can be varied.
Various cases can be directly observed from Fig. \ref{fig: interscept} by comparing the linear function on the
left-hand side and $-F(u)$ on the right-hand side. For two particular values of $P$, the linear function becomes tangent to $-F(u)$. They correspond to two bifurcation points
of saddle-node type. In spite of the simplicity of the one-dimensional lattice model, the bifurcation seems to be quite generic. In fact,  the same type of bifurcations have been observed in two and three-dimensional lattice models \cite{XLi2013A}, where a sequence of saddle-node bifurcations were observed.

In what follows, we will turn to the QC approximation models, and investigate how the bifurcation diagram is influenced by the QC approximation.
\section{Crack-tip Solutions and Bifurcation Curves for the Multiscale Models}\label{sec: multi}

We analyze three QC methods applied to the above lattice model. The calculation will be carried out as explicitly as possible, with the goal of not overestimating or underestimating the error.
Due to the discrete nature of the model, the calculation is quite lengthy. We will only show the
 full details for the first model, and keep the procedure brief for the other two models.
\subsection{The quasicontinuum method without force correction}
The original
QC method~\cite{TadmorOrtizPhillips:1996} relies on an energy summation rule. In the cracked region, the total energy can be written as a sum of the site energy, i.e.,
\( E= \sum_i V_i\) with
\[
V_i =\dfrac{\ka_1}{4}\Lr{(u_{j+1}-u_j)^2+(u_{j-1}-u_j)^2}+ \dfrac{\ka_2}{4}\Lr{(u_{j+2}-u_j)^2 + (u_{j-2}-u_j)^2}.
\]
Moreover, for the atoms at and behind the crack-tip, an energy functional
for the vertical bonds should be included in the total energy.

The QC method introduces a {\em local} region, where the displacement field is represented on a finite element mesh, and within each element, the energy is approximated by Cauchy-Born (CB) rule~\cite{BornHuang:1954}. To separate out the issue of interpolation and quadrature error, we assume that the mesh node coincides with the atom position. In this case, the approximating energy takes the form of \(E_{\rm QC} = \sum_i E_i\), where the summation is over all the atom sites. We assume that the local region includes atoms \(j<m-1\), and  the approximating energy is given by
\[
E_j= \dfrac{\ov{\ka}}{2}\bigl[(u_{j+1}-u_j)^2+(u_{j}-u_{j-1})^2\bigr]\qquad
\text{for\quad}j<m-1
\]
as a result of the CB approximation. In addition, we have
\[
E_j=V_j\qquad\text{for\quad}m<j<n.
\]
At the interface, the energy functions are given by
\begin{align*}
E_{m-1}&=\dfrac{\ka_1}{4}\Lr{(u_m - u_{m-1})^2+(u_{m-1} - u_{m-2})^2}+\ka_2\Lr{(u_m - u_{m-1})^2 + (u_{m-1} - u_{m-2})^2},\\
E_m &=\dfrac{\ka_1}{4}
\Lr{(u_m - u_{m-1})^2+(u_{m+1} - u_{m})^2}+\dfrac{\ka_2}{4}
\Lr{(u_{m+2} - u_{m})^2+4(u_m - u_{m-1})^2}.
\end{align*}
We have the following system of equilibrium equations:
\[
\left\{
\begin{aligned}
\ov{\ka}(u_1-u_0) + P&=0,\\
\ov{\ka}(u_{j+1}-2u_j+u_{j-1})&=0,\qquad 2\le j \le m-2,\\
 L^A u_j&=0,\qquad m-1\le  j \le n-1,\\
L^Au_n+F(u_n)&=0,\\
L^Au_j+2\ka_3 u_j&=0,\qquad j\ge n+1.
\end{aligned}\right.
\]
Around the interface, we have the following coupling equations:
\begin{equation}\label{eq:qccouple0}
\left\{
\begin{aligned}
\ov{\ka}u_{m-2} - \Lr{2\ka_1+17\ka_2/2}u_{m-1} +\ov{\ka}u_{m}+ \dfrac{\ka_2}2u_{m+1}&=0,\\
 \ov{\ka}u_{m-1} - (2\ka_1+5\ka_2)u_m+ \ka_1 u_{m+1} + \ka_2 u_{m+2}&=0,\\
\dfrac{\ka_2}{2}u_{m-1} +\ka_1u_{m}-\Lr{2\ka_1+3\ka_2/2}u_{m+1}+ \ka_1 u_{m+2}+ \ka_2 u_{m+3}&=0.
\end{aligned}\right.
\end{equation}

In the local region, 
\[
  u_j = C_0 + C_1 j, \quad j=0, 1, \cdots, m-1
\]
for certain constants $C_0$ and $C_1$. If we impose the traction boundary condition, then the solution takes a simpler form as
\begin{equation}\label{eq:linearcb}
u_j=C_0+\dfrac{P}{\ov{\ka}}(m-j-1).
\end{equation}

Adding up all the equations in~\eqref{eq:qccouple0}, we obtain
\begin{equation}\label{eq:qccouple1}
(\ka_1+\ka_2)(u_{m+2}-u_{m+1})+\ka_2(u_{m+3}-u_m)=-\ov{\ka}(u_{m-2}-u_{m-1})
=-P,
\end{equation}
where we have used~\eqref{eq:linearcb} in the last step.

We substitute~\eqref{eq:linearcb} into the first
two equations of~\eqref{eq:qccouple0} and obtain
\begin{align*}
(\ka_1 +9\ka_2/2)(u_m- u_{m-1})-\dfrac{\ka_2}{2}(u_{m+1}-u_m)&=-P,\\
 -\ov{\ka}(u_m-u_{m-1})+ \ka_1(u_{m+1}-u_m) + \ka_2(u_{m+2}- u_m)&=0.
\end{align*}
Denote \(\gamma=\ov{\ka}/[\ov{\ka}+\ka_2/2]\). We eliminate $u_m-u_{m-1}$ from the above two equations and obtain the following linear system.
\begin{equation}\label{eq:couple}
\left\{\begin{aligned}
\Bigl[ \ka_1 + \dfrac{\gamma}{2}\ka_2\Bigr] (u_{m+1}-u_m) + \ka_2(u_{m+2}- u_m)
 &=-\gamma P,\\
(\ka_1 + \ka_2)(u_{m+2}-u_{m+1})+ \ka_2 (u_{m+3}-u_{m+1})&=-P.
 \end{aligned}\right.
\end{equation}

To proceed, we express the solution in the atomistic region before the crack-tip in the same form as in the previous section, for \(j=m-3,\cdots,n+1\),
\begin{equation}\label{eq:atomansatz}
u_j=a+bj+c\cosh[j\del]+d\sinh[j\del].
\end{equation}
We substitute the above ansatz into~\eqref{eq:couple} and obtain
\begin{equation}\label{eq:qccouple}
(c,d)\mc{K}_m=(P,b)\mc{Q},
\end{equation}
where $\mc{Q}=\{q_{ij}\}_{i,j=1}^2$ is a $2$ by $2$ matrix given by
\[
\mc{Q}{:}=\dfrac1{4\ka_2}\begin{pmatrix}
\dfrac{4-3\ga}{\cosh\del-1}
&\dfrac{3\ga}{\sinh\del}\\[.1in]
\dfrac{(2-\ga)\ov{\ka}}{\cosh\del-1}&\dfrac{(\ga+2)\ov{\ka}-4\ka_2}{\sinh\del}
\end{pmatrix}.
\]
The details are postponed to Appendix~\ref{app:qc}.

Using the fact that \(\mc{K}_n\mc{K}_m^{-1}=\mc{K}_{n-m}\), we get
\begin{align*}
\mc{F}_{n,\al}(\del)c+\mc{G}_{n,\al}(\del)d&=(c,d)\mc{K}_n(A_\al,B_\al)^T
=(P,b)\mc{Q}\mc{K}_m^{-1}\mc{K}_n
(A_\al,B_\al)^T\\
&=(P,b)\mc{Q}\mc{K}_{n-m}(A_\al,B_\al)^T.
\end{align*}
Using~\eqref{eq:matchcoef}, we
represent $b$ in terms of $u_n$ and $P$ as
\[
b=\dfrac{\al+\beta-1}{q_{12}\mc{F}_{n-m,\al}
+q_{22}\mc{G}_{n-m,\al}+1+\al}u_n-\dfrac{q_{11}\mc{F}_{n-m,\al}
+q_{21}\mc{G}_{n-m,\al}}
{q_{12}\mc{F}_{n-m,\al}+
q_{22}\mc{G}_{n-m,\al}+1+\al}P.
\]
A direct calculation gives
\begin{align}
\mc{F}_{n,1-\beta}(\del)c+\mc{G}_{n,1-\beta}(\del)d
&=\Lr{q_{12}\mc{F}_{n-m,1-\beta}
+q_{22}\mc{G}_{n-m,1-\beta}}b\nn\\
&\quad+\Lr{q_{11}\mc{F}_{n-m,1-\beta}+
q_{21}\mc{G}_{n-m,1-\beta}}P.\label{eq:basrela}
\end{align}

Now we find the effective equation for $u_n$:
\begin{align*}
F(u)+\ka^{\rm qc} u+\eta^{\rm qc} P=0
\end{align*}
with
\begin{align*}
\ka^{\rm qc}&=(\al+\beta-1)\Lr{\ka_1+\ka_2(1+\beta)
+\ka_2\dfrac{
q_{12}\mc{F}_{n-m,1-\beta}+
q_{22}\mc{G}_{n-m,1-\beta}}
{q_{12}\mc{F}_{n-m,\al}+q_{22}\mc{G}_{n-m,\al}+1+\al}}\\
&\quad-(\al+\beta-1)\dfrac{
\ov{\ka}+(\beta-2)\ka_2}
{q_{12}\mc{F}_{n-m,\al}+q_{22}\mc{G}_{n-m,\al}+1+\al},
\end{align*}
and
\begin{align*}
\eta^{\rm qc}&=\ka_2\Lr{q_{11}\mc{F}_{n-m,1-\beta}
+q_{21}\mc{G}_{n-m,1-\beta}}\\
&\quad-\ka_2\dfrac{\Lr{q_{11}\mc{F}_{n-m,\al}
+q_{21}\mc{G}_{n-m,\al}}\Lr{q_{12}\mc{F}_{n-m,1-\beta}
+q_{22}\mc{G}_{1-\beta}}}
{q_{12}\mc{F}_{n-m,\al}+q_{22}\mc{G}_{n-m,\al}+1+\al}\\
&\quad+\dfrac{\Lr{q_{11}\mc{F}_{n-m,\al}
+q_{21}\mc{G}_{n-m,\al}}
\bigl(\ov{\ka}+(\beta-2)\ka_2\bigr)}
{q_{12}\mc{F}_{n-m,\al}
+q_{22}\mc{G}_{n-m,\al}+1+\al}.
\end{align*}

Let $n-m\to\infty$, we obtain $\ka^{\text{qc}}\to\ka_0$ with the expansion
\[
\ka^{\text{qc}}=\ka_0^{\text{qc}}+\dfrac{2(\al+\beta-1)(\al+\beta-1-A_\al-B_\al)
\ov{\ka}}
{(q_{12}+q_{22})(A_\al+B_\al)^2}z_0^{n-m}+\mc{O}(z_0^{2(n-m)}).
\]
Proceeding along the same line that leads to~\eqref{eq:crisscross}, we obtain
\begin{align*}
&\Lr{q_{11}\mc{F}_{n-m,1-\beta}
+q_{21}\mc{G}_{n-m,1-\beta}}
\Lr{q_{12}\mc{F}_{n-m,\al}+q_{22}\mc{G}_{n-m,\al}}\nn\\
&\quad-\Lr{q_{11}\mc{F}_{n-m,\al}+q_{21}\mc{G}_{n-m,\al}}
\Lr{q_{12}\mc{F}_{n-m,1-\beta}+q_{22}\mc{G}_{n-m,1-\beta}}\nn\\
&=\det\biggl[\mc{Q}\begin{pmatrix}
\mc{F}_{n-m,1-\beta}&\mc{G}_{n-m,1-\beta}\\
\mc{F}_{n-m,\al}&\mc{G}_{n-m,\al}
\end{pmatrix}\biggr]\nn\\
&=\det\mc{Q}\det\begin{pmatrix}
\mc{F}_{n-m,1-\beta}&\mc{G}_{n-m,1-\beta}\\
\mc{F}_{n-m,\al}&\mc{G}_{n-m,\al}
\end{pmatrix}\nn\\
&=-(\al+\beta-1)\Lr{2(1-\ga)\ov{\ka}-(4-3\ga)\ka_2}/(2\ka_2^2).
\end{align*}

Using the above identity, we write $\eta^{\rm qc}$ as
\begin{align*}
\eta^{\rm qc}&=\dfrac{q_{11}\mc{F}_{n-m,\al}
+q_{21}\mc{G}_{n-m,\al}}
{q_{12}\mc{F}_{n-m,\al}
+q_{22}\mc{G}_{n-m,\al}+1+\al}\bigl(\ov{\ka}+(\beta-2)\ka_2\bigr)\\
&\quad+
\dfrac{q_{11}\mc{F}_{n-m,1-\beta}
+q_{21}\mc{G}_{n-m,1-\beta}}
{q_{12}\mc{F}_{n-m,\al}
+q_{22}\mc{G}_{n-m,\al}+1+\al}(1+\al)\ka_2\\
&\quad-\dfrac{(\al+\beta-1)\Lr{2(1-\ga)\ov{\ka}-(4-3\ga)\ka_2}}
{2\Lr{q_{12}\mc{F}_{n-m,\al}
+q_{22}\mc{G}_{n-m,\al}+1+\al}\ka_2},
\end{align*}
which can be reshaped into
\begin{align*}
\eta^{\rm qc}&=\dfrac{q_{11}\mc{F}_{n-m,\al}
+q_{21}\mc{G}_{n-m,\al}}
{q_{12}\mc{F}_{n-m,\al}+q_{22}\mc{G}_{n-m,\al}+1+\al}\ov{\ka}\\
&\quad+\ka_2
\dfrac{q_{11}\Lr{(\beta-2)\mc{F}_{n-m,\al}+(1+\al)\mc{F}_{n-m,1-\beta}}
+q_{21}\Lr{(\beta-2)\mc{G}_{n-m,\al}+(1+\al)\mc{G}_{n-m,1-\beta}}}
{q_{12}\mc{F}_{n-m,\al}+q_{22}\mc{G}_{n-m,\al}+1+\al}\\
&\quad-\dfrac{(\al+\beta-1)\Lr{(1-\ga)\ov{\ka}/\ka_2-(4-3\ga)}}
{q_{12}\mc{F}_{n-m,\al}
+q_{22}\mc{G}_{n-m,\al}+1+\al}.
\end{align*}
Next we proceed with the same procedure that leads to~\eqref{eq:idencoeff}, and obtain
\begin{align*}
\eta^{\rm qc}&=\dfrac{q_{11}\mc{F}_{n-m,\al}
+q_{21}\mc{G}_{n-m,\al}}
{q_{12}\mc{F}_{n-m,\al}+q_{22}\mc{G}_{n-m,\al}+1+\al}\ov{\ka}\\
&\quad-(\al+\beta-1)\ov{\ka}
\dfrac{q_{11}\cosh[(n-m)\del]+q_{21}\sinh[(n-m)\del]}
{q_{12}\mc{F}_{n-m,\al}+q_{22}\mc{G}_{n-m,\al}+1+\al}\\
&\quad-\dfrac{(\al+\beta-1)\Lr{(1-\ga)\ov{\ka}/\ka_2-(2-3\ga/2)}}
{q_{12}\mc{F}_{n-m,\al}
+q-{22}\mc{G}_{n-m,\al}+1+\al}.
\end{align*}
Let $n-m\to\infty$, we obtain $\eta^{\text{qc}}\to\eta_0^{\text{qc}}$ with
\begin{align*}
\eta_0^{\text{qc}}&=\Lr{1-\dfrac{\al+\beta-1}{A_\al+B_\al}}\dfrac{q_{11}+q_{21}}{q_{12}+q_{22}}\ov{\ka}\\
&=\Lr{1-\dfrac{\al+\beta-1}{A_\al+B_\al}}\dfrac{4+3\ga(\tanh[\del/2]-1)}
{2-\ga+(\ga+2-4\ka_2/\ov{\ka})\tanh[\del/2]}.
\end{align*}
A direct calculation gives
\[
\eta_0-\eta_0^{\text{qc}}=\Lr{1-\dfrac{\al+\beta-1}{A_\al+B_\al}}
\Lr{1-\dfrac{q_{11}+q_{21}}{q_{12}+q_{22}}\ov{\ka}}.
\]
It is clear that the difference remains finite as $n-m\to\infty$. Therefore,
$\eta_0^\mathrm{qc}$ does not
coincide with \(\eta_0\).
\medskip

\begin{figure}[thbp]
\begin{center}
\includegraphics[scale=0.5]{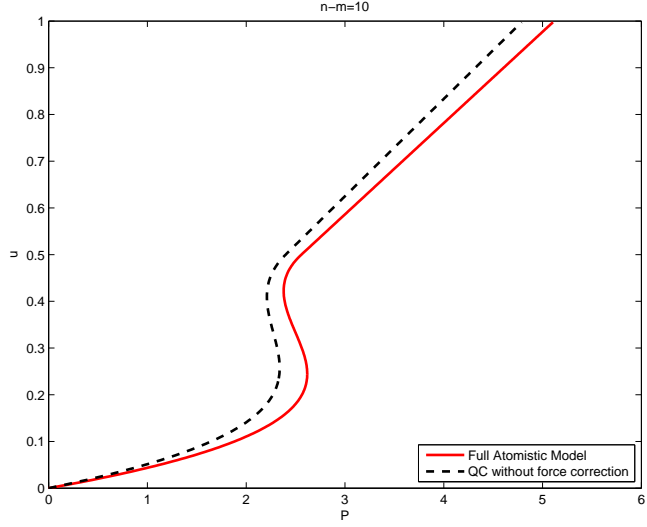}
\caption{The bifurcation diagram for the full atomistic model and QC without force corrections. The middle branch contains unstable equilibrium while the other two branches are stable. }
\label{fig: bif_qc_10}
\end{center}
\end{figure}

As an example of comparison, we plot the bifurcation diagram for both models in Fig. \ref{fig: bif_qc_10}. We chose \(\ka_1=4\),
\(\ka_2=0.4\), \(\ka_3=20\), and \(u_{\text{cut}}=0.5\). Clearly, the diagram consists of two saddle-node bifurcation points. Although QC predicts a similar bifurcation behavior, the error in the bifurcation curves is significant.
\subsection{Quasinonlocal Quasicontinuum method}
The quasinonlocal quasicontinuum method (QQC)~\cite{Shimokawa:2004} approximates the energy as follows. For \(j<m-1\), the site energy is
\[
E_j=\dfrac{\ov{\ka}}{2}\bigl((u_{j+1}-u_j)^2+(u_{j}-u_{j-1})^2\bigr),
\]
and for \(j>m\), we set
\(E_j = V_j\), and at the interface, i.e., for $j=m-1,m$,
\[
  E_j=\dfrac{\ka_2}{4} (u_{j+2}-u_j)^2 + \dfrac{\ka_1}{4}\bigl((u_{j+1} - u_j)^2
  +(u_j - u_{j-1})^2\bigr) +\ka_2(u_j-u_{j-1})^2.
\]

The resulting force balance equations take the following form:
\[
\left\{
\begin{aligned}
\ov{\ka}(u_1 - u_0)  + P &=0,\\
\ov{\ka}(u_{j-1}-u_j) + \ov{\ka}(u_{j+1}-u_j)&=0,\qquad
1 \le j \le m-2\\
L^Au_j&=0,\quad m < j < n,
\end{aligned}\right.
\]
and around the interface,
\begin{equation}\label{eq:qqcinter}
\left\{
\begin{aligned}
\ov{\ka}(u_{m-2}-u_{m-1}) + (\ka_1+2\ka_2)(u_m - u_{m-1}) + \ka_2(u_{m+1}- u_{m-1}) &= 0,\\
(\ka_1 + 2\ka_2)(u_{m-1}-u_m)+\ka_1(u_{m+1}-u_m) + \ka_2(u_{m+2}- u_m)&=0,
\end{aligned}\right.
\end{equation}
Using a similar procedure that leads to~\eqref{eq:couple}, we eliminate \(u_{m-1}-u_{m-2}\) from~\eqref{eq:qqcinter}
and obtain
\begin{equation}\label{eq: qqc'}
\left\{
 \begin{aligned}
(\ka_1+2\ka_2)(u_m - u_{m-1}) + \ka_2(u_{m+1}-u_{m-1}) &= -P, \\
 -(\ka_1+2\ka_2)(u_m - u_{m-1}) + \ka_1(u_{m+1}-u_m) + \ka_2(u_{m+2}-u_{m})&=0.
 \end{aligned}\right.
\end{equation}
Substituting the general expression of $u_j$ into the above two
equations, we obtain
\begin{equation}\label{eq:qqc}
(c,d)\mc{K}_{m-1}=(P/\ov{\ka}+b,0).
\end{equation}
We leave the details for deriving the above equation to the
Appendix.
Solving the above equation, we obtain
\(b=-P/\ov{\ka}\),
which together with~\eqref{eq: matching-tip} yields
\[
\mc{F}_n(\del)c+\mc{G}_n(\del)d=(\al+\beta-1)u_n+(1+\al)P/\bar{\ka}.
\]
This equation, together with~\eqref{eq:qqc}$_2$, gives
\[
\left\{
\begin{aligned}
c&=-\dfrac{\sinh[(m-1)\del]}{\mc{G}_{n-m+1,\al}(\del)}
\bigl((\al+\beta-1)u_n+(1+\al)P/\bar{\ka}\bigr),\\
d&=\dfrac{\cosh[(m-1)\del]}{\mc{G}_{n-m+1,\al}(\del)}
\bigl((\al+\beta-1)u_n+(1+\al)P/\bar{\ka}\bigr).
\end{aligned}\right.
\]

Substituting the expressions of $c$ and $d$ into~\eqref{eq:matchcoef}, we obtain
\[
F(u_n)+\ka^{\text{qqc}} u_n+\eta^{\text{qqc}}P=0
\]
with
\begin{align*}
\ka^{\text{qqc}}&=(\al+\beta-1)\Lr{\ka_1+\ka_2(\beta+1)+\ka_2
\dfrac{\mc{G}_{n-m+1,1-\beta}(\del)}{\mc{G}_{n-m+1,\al}(\del)}},\\
\eta^{\text{qqc}}&=1-(\al+\beta-1)\dfrac{\sinh[(n-m+1)\del]}
{\mc{G}_{n-m+1,\al}(\del)}.
\end{align*}

We expand these two parameters and get
\begin{equation}\label{eq: kappa-exp-qqc}
\ka^{\text{qqc}}=\ka_0 -\dfrac{2\ov{\ka} (\alpha+\beta-1)\sinh\del}{(A_\al+B_\al)^2}z_0^{2n-2m+2}
+\mc{O}(z_0^{4n-4m+4}),
\end{equation}
and
\begin{equation}\label{eq: eta-exp-qqc}
\eta^{\text{qqc}}=\eta_0 + \dfrac{2(\alpha+\beta-1)A_\al}{(A_\al+B_\al)^2}z_0^{2n-2m+2}
+\mc{O}(z_0^{4n-4m+4}).
\end{equation}

\begin{figure}[htbp]
\begin{center}
\includegraphics[scale=0.5]{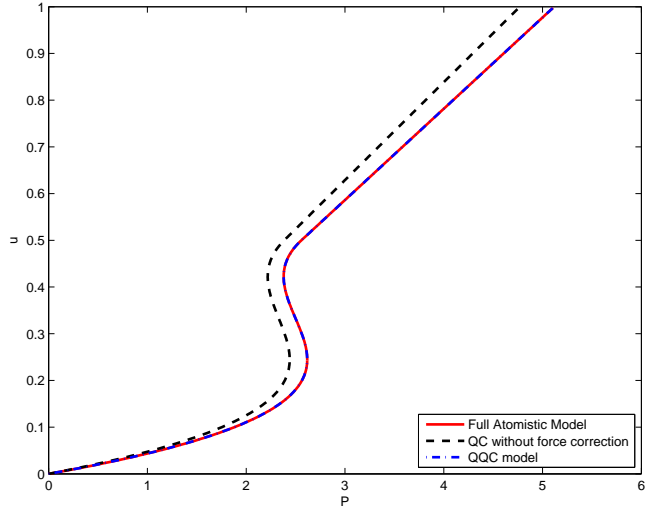}
\caption{The bifurcation diagram for the QC and QQC method. }
\label{fig: bif_qc_qqc}
\end{center}
\end{figure}

Let $n-m\to\infty$, we obtain
\(\ka\to\ka_0\) and \(\eta\to\eta_0.\)
The limits $\ka_0$ and $\eta_0$ are the same with those of the atomistic model.
Namely, QQC gives the correct bifurcation graph in the limit.
This can be confirmed from a comparison of the bifurcation diagram, as illustrated in Fig. \ref{fig: bif_qc_qqc}, where excellent agreement is observed. To reach this asymptotic regime, the crack-tip has to be sufficiently far away from the atomistic/continuum interface.
\subsection{A force-based method}
A force-based quasicontinuum (FQC) method differs from the previous methods
lies in the fact that there does not exist an associated energy. One simple approach to construct a FQC method is to keep the equations~\eqref{eq: cracked} through
\eqref{eq: un-cracked} in the atomistic region, while the force balance equations computed from the Cauchy-Born rule are used in the continuum region. The resulting equilibrium equations read as
\[
\left\{
\begin{aligned}
\ov{\ka}(u_1-u_0) + P&=0, \\
\ov{\ka}(u_{j+1}-2u_j+u_{j-1}) &=0, \qquad 2\le j \le m,\\
L^A u_j&=0,\qquad m+1\le j < n \\
L^Au_n+F(u_n)&=0,\\
L^Au_j+2\ka_3u_j&=0,\qquad j\ge n+1.
\end{aligned}\right.
\]

In this case, we still have
\begin{equation}\label{eq: fb'}
u_m= u_{m+1} + P/\ov{\ka},\quad u_m= u_{m+2} +2P/\ov{\ka}.
\end{equation}

Substituting the expression of $u_j$ into the above equation, we obtain
\begin{equation}\label{eq:fqccouple}
(c,d)\mc{K}_{m+1}=-(\coth\del,1)(P/\ov{\ka}+b).
\end{equation}
Solving the above equations, we get
\begin{equation}\label{eq:cd}
c=\dfrac{\sinh[m\del]}{\sinh\del}(P/\bar{\ka}+b),\qquad
d=-\dfrac{\cosh[m\del]}{\sinh\del}(P/\bar{\ka}+b).
\end{equation}
Substituting the expressions of $c$ and $d$ into~\eqref{eq:matchcoef}, we obtain
\[
b=\dfrac{\mc{G}_{n-m,\al}(\del)}{(1+\al)\sinh\del-\mc{G}_{n-m,\al}(\del)}
\dfrac{P}{\ov{\ka}}+\dfrac{(\al+\beta-1)\sinh\del}
{(1+\al)\sinh\del-\mc{G}_{n-m,\al}(\del)}u_n.
\]
Finally, we substitute the expressions of $b,c$ and $d$ into~\eqref{eq:match-expa}
and obtain
\[
F(u_n)+\ka^{\text{fqc}} u_n+\eta^{\text{fqc}} P=0
\]
with
\begin{align*}
\ka^{\text{fqc}}&=(\al+\beta-1)\Lr{\ka_1+\ka_2(\beta+1)
+\ka_2\dfrac{\mc{G}_{n-m,1-\beta}(\del)}
{\mc{G}_{n-m,\al}(\del)-(1+\al)\sinh\del}}\\
&\quad+(\al+\beta-1)\dfrac{(\ov{\ka}+(\beta-2)\ka_2)\sinh\del}
{\mc{G}_{n-m,\al}(\del)-(1+\al)\sinh\del},
\end{align*}
and
\begin{align*}
\eta^{\text{fqc}}&=\Lr{1+\dfrac{(\beta-2)\ka_2}{\ov{\ka}}}
\Lr{1+\dfrac{(1+\al)\sinh\del}{\mc{G}_{n-m,\al}(\del)-(1+\al)\sinh\del}}\\
&\quad+\dfrac{(1+\al)\ka_2}{\ov{\ka}}\dfrac{\mc{G}_{n-m,1-\beta}(\del)}
{\mc{G}_{n-m,\al}(\del)(1+\al)\sinh\del}.
\end{align*}

We expand \(\ka^{\text{fqc}}\) as follows.
\begin{equation}\label{eq: kappa-exp-fb}
\ka^{\text{fqc}}=\ka_0 + 2(\al+\beta-1)\ov{\ka}\dfrac{A_\al+B_\al-(\al+\beta-1)\sinh\del}
{(A_\al+B_\al)^2}z_0^{n-m} + \mc{O}(z_0^{2n-2m}).
\end{equation}

We write $\eta^{\text{fqc}}$ as
\begin{align*}
\eta^{\text{fqc}}=1-\dfrac{(\al+\beta-1)\sinh\del\sinh[(n-m)\del]}
{\mc{G}_{n-m,\al}(\del)-(1+\al)\sinh\del}
+\dfrac{(1+\al)\sinh\del}{\mc{G}_{n-m,\al}(\del)-(1+\al)\sinh\del}.
\end{align*}
Hence we have
\begin{equation}\label{eq: eta-exp-fb}
\eta^{\text{fqc}}=\eta_0 + \dfrac{(3+\al-\beta)\sinh\del}{A_\al+B_\al}
z_0^{n-m} + \mc{O}(z_0^{2n-2m}).
\end{equation}

It is clear that
\(\ka\to\ka_0\) and \(\eta \to \eta_0\) as $n-m\to\infty$.
\subsection{A comparison test}
As an example, we continue from the first numerical test and set \(m=100\) and \(n=104\). We computed the coefficients and listed the results in the following table.
It is clear that QQC gives the best results, while the QC method gives a good approximation of \(\kappa\) but the error in \(\eta\) is large.
%
\begin{table}[h]\label{tab: test1}
  \begin{center}
    \begin{tabular}{|c|c|c|}
        \hline
         & $\kappa$ & $\eta$ \\
            \hline
        Exact & -4.782062040603841 &  0.934371338155818 \\
        QC &   -4.782048350329799  & 1.002417909367481 \\
	QQC&  -4.782060913687936  & 0.934371132296173 \\
         FQC & -4.782243406077938 &  0.934404469139225 \\
        \hline
    \end{tabular}
    \end{center}
    \caption{The parameters for the effective equations when \(m=100\) and \(n=104\)}.
\end{table}

In the second test, we set \(m=96\) while keeping $n=104$. This corresponds to enlarging the atomistic region. The results are collected in
Table~\ref{tab: test2}. We notice that for the QC method, the parameter $\kappa$ has further approached to the exact value. However, the error in $\eta$ still remains finite,
which confirms our analysis. On the other hand, the error for QQC and FQC have
been greatly improved.
\begin{table}[h]\label{tab: test2}
  \begin{center}
    \begin{tabular}{|c|c|c|}
        \hline
         & $\kappa$ & $\eta$ \\
            \hline
        Exact & -4.782062040603841 &  0.934371338155818 \\
        QC &   -4.782062040081748  & 1.002420436153853  \\
	QQC&  -4.782062040560865  & 0.934371338147967\\
         FQC&  -4.782062047519995 &  0.934371339419228 \\
        \hline
    \end{tabular}
    \end{center}
    \caption{The parameters for the effective equations when \(m=96\) and \(n=104\)}.
\end{table}
\section{Analysis of the Bifurcation Curves}\label{sec: ode}
Now we are ready to  estimate the overall error, and we will focus on
the error in $u_n$ because the displacement of any other atoms can be expressed as a linear function of $u_n$ as shown in the previous section. The
error for $u_n$ is the error committed by solving the effective equations with the approximating parameters $\kappa$ and $\eta$. It follows from Fig.~\ref{fig: bif_qc_10} and Fig.~\ref{fig: bif_qc_qqc} that there might be multiple solutions of the effective equation
with a given load $P$. In addition,
there are two points when the derivative with respect to $P$ is infinite. These two points are exactly the bifurcation points. Therefore, it is difficult to compare $u_n$ with its approximations directly for the same loading parameter. In fact, we expect the error of $u_n$ to be quite large near the bifurcation points.

Instead of a direct comparison, we propose a different approach, which is motivated by continuation methods for solving bifurcation problems~\cite{RhBu83,Seydel10}. More specifically, we will compare the bifurcation curves as a whole. For this purpose, we
parameterize the bifurcation curve on the $P-u$ plane using arc length, which is denoted by $s$. Compared to the parameterization using the load parameter, the new representation is not multi-valued. First we set the initial point of the curve
to \( (0,0)\), which clearly satisfies the effective equation for any choice of the parameters $\kappa$ and $\eta$. Next we represent a point on the curve by $(P(s), u(s))$. To trace out the curve, one compute the tangent vector
\[
 \tau(s)= \big(f_1(u(s); \ka, \eta),f_2(u(s); \ka, \eta)\big),
\]
where
\[
f_1\big(u(s); \ka, \eta\big)= -\dfrac{\eta}{\sqrt{(F'(u(s))+\ka)^2+\eta^2}}, \quad
f_2\big(u(s); \ka, \eta\big) =\dfrac{F'(u(s))+\ka}{\sqrt{(F'(u(s))+\ka)^2+\eta^2}}.
\]
This can be easily obtained by differentiation the effective equation with respect to the arc length.
Following the curve with $s$ as the independent variable, we obtain the following ODEs that describe the bifurcation curve~\cite{RhBu83}:
\[
\left\{ \begin{aligned}
   \frac{d}{ds}u(s)&= f_1\big(u(s); \ka, \eta\big),\\
   \frac{d}{ds}P(s)&= f_2\big(u(s); \ka, \eta\big),\\
   \quad u(0)&=0, \; P(0)=0.
 \end{aligned}\right.
\]

As we have shown in the previous sections, a multiscale method typically gives an effective equation for $u_n$ that is of the same form as the exact equation, but with the approximate parameters $\ka$ and $\eta$.
We denote the approximated values as $\widehat{\ka}$ and $\widehat{\eta}$, and the corresponding bifurcation curve as $\big(\widehat{P}(s), \widehat{u}(s)\big)$, respectively. We can describe
 the bifurcation curve by the following ODEs:
 \[
 \left\{ \begin{aligned}
   \frac{d}{ds}\widehat{u}(s)&= f_1\big(\widehat{u}(s); \widehat{\ka}, \widehat{\eta}\big),\\
   \frac{d}{ds}\widehat{P}(s)&= f_2\big(\widehat{u}(s); \widehat{\ka}, \widehat{\eta}\big),\\
   \quad \widehat{u}(0)&=0, \; \widehat{P}(0)=0.
 \end{aligned}\right.
\]
In this way, the problem has been reduced to a perturbation problem with varying parameters.
Standard theory for ODEs states that the solution is continuously dependent on the parameters~\cite{Perko}, provided that the functions $f_1$ and $f_2$ are Lipschitz continuous. This can be explicitly stated as follows.
For any $s\in [0, S]$,
\[
 \abs{u(s) - \widehat{u}(s)}+ \abs{P(s)-\widehat{P}(s)}\le L\Lr{\abs{\ka - \widehat{\ka}}+\abs{\eta - \widehat{\eta}}}e^{LS},
\]
were $L$ is the Lipschitz constant of $f_1$ and $f_2$. In particular, the error in $u_n$ will depend continuously on the parameters $\kappa$ and $\eta$, and for the QQC and the force-based methods, this error should be exponentially small.   More importantly, this estimate is not restricted to a local minimum.
\section{Summary}
In this paper,
We have evaluated the error of three QC methods by comparing the bifurcation diagram. In particular, it has been found that the original QC method, with the notorious problem of ghost forces, exhibits large error in
predicting the bifurcation curve. The quasi-nonlocal QC method and a force-based method, on the
other hand, are quite accurate in this aspect. This suggests that ghost forces are responsible for the large error.  Using the parametrization with arc length, we have been able to obtain {\it quantitative} estimate for the approximation of the bifurcation curves.

One remaining issue is estimating the error in the continuum region. For the current problem, once
$u_n$ is obtained from the bifurcation diagram, the rest of the degrees of freedom are uniquely determined. This makes it possible to interpret the error in the continuum region. In the context
of bifurcation theory, the effective equation~\eqref{eq: eff0} describes a center manifold, where the transition occurs. The remaining degrees of freedom lie in the stable manifold, and  standard methods in numerical analysis
may apply. This issue for more general problems will be addressed in futures works.
\begin{appendix}
\section{Derivation of Equations~\eqref{eq:qccouple}}\label{app:qc}
We first introduce a shorthand notation. For $a,x\in\mb{R}$, denote
\[
s_a(x)=\sinh[ax],\quad c_a(x)=\cosh[ax].
\]
\vskip .5cm
\noindent{\em Proof of ~\eqref{eq:rela1} and~\eqref{eq:rela2}.\;}
The identity~\eqref{eq:rela1} is equivalent to
\[
\ka_2\Lr{c_{k+1}(\del)-c_{k-2}(\del)}+(\ka_1+\ka_2)\Lr{c_k(\del)-c_{k-1}(\del)}=0.
\]
The left-hand side of the above equation can be written into
\[
2s_{k-1/2}(\del)\Lr{\ka_2s_{3/2}(\del)+(\ka_1+\ka_2)s_{1/2}(\del)}.
\]
Using~\eqref{eq:pararela1}, we obtain
\[
\ka_2s_{3/2}(\del)+(\ka_1+\ka_2)s_{1/2}(\del)=s_{1/2}(\del)\Lr{\ka_2(2c_1(\del)+1)+
\ka_1+\ka_2}=0.
\]
This completes the proof for~\eqref{eq:rela1}.

We omit the proof for~\eqref{eq:rela2} since it is the same.$\blacksquare$

First we substitute~\eqref{eq:atomansatz} into~\eqref{eq:couple}$_2$ and obtain
\begin{align*}
&\Bigl\{(\ka_1+\ka_2)\Lr{c_{m+2}(\del)-c_{m+1}(\del)}
+\ka_2\Lr{c_{m+3}(\del)-c_{m+1}(\del)}\Bigr\}c\\
&+\Bigl\{(\ka_1+\ka_2)\Lr{s_{m+2}(\del)-s_{m+1}(\del)}
+\ka_2\Lr{s_{m+3}(\del)-s_{m+1}(\del)}\Bigr\}d\\
&=-P-(\ka_1+3\ka_2)b.
\end{align*}
Using~\eqref{eq:rela1} and~\eqref{eq:rela1}with $k=m+2$ to simplify the coefficients
for $c$ and $d$, respectively, we obtain
a simplified form of~\eqref{eq:couple}$_2$ as
\begin{equation}\label{eq:couple2}
\bigl(c_{m+1}(\del)-c_m(\del)\bigr)c+\bigl(s_{m+1}(\del)-s_m(\del)\bigr)d
=\dfrac{P}{\ka_2}+\dfrac{\ka_1+3\ka_2}{\ka_2}b.
\end{equation}
Next we substitute~\eqref{eq:atomansatz} into~\eqref{eq:couple}$_1$ and obtain
\begin{equation}\label{eq:couple1}
\begin{aligned}
&\Bigl\{(\ka_1+\ga\ka_2/2)\Lr{c_{m+1}(\del)-c_m(\del)}
+\ka_2\Lr{c_{m+2}(\del)-c_m(\del)}\Bigr\}c\\
&+\Bigl\{(\ka_1+\ga\ka_2/2)\Lr{s_{m+1}(\del)-s_m(\del)}
+\ka_2\Lr{s_{m+2}(\del)-s_m(\del)}\Bigr\}d\\
&+\Lr{\ka_1+(\ga/2+2)\ka_2}b=-\ga P.
\end{aligned}
\end{equation}

A direct calculation yields
\begin{align*}
\ka_1\Lr{c_{m+1}(\del)-c_m(\del)}+\ka_2\Lr{c_{m+2}(\del)-c_m(\del)}
&=-2\ka_2s_1(\del)s_m(\del),\\
\ka_1\Lr{s_{m+1}(\del)-s_m(\del)}+\ka_2\Lr{s_{m+2}(\del)-s_m(\del)}
&=-2\ka_2s_1(\del)c_m(\del).
\end{align*}
Using the above two equations, we  may write the equation~\eqref{eq:couple1} into
\[
\begin{aligned}
&\Lr{\dfrac\ga{2}\Lr{c_{m+1}(\del)-c_m(\del)}-2s_1(\del)s_m(\del)}c
+\Lr{\dfrac\ga{2}\Lr{s_{m+1}(\del)-s_m(\del)}-2s_1(\del)c_m(\del)}d\\
&=-\dfrac{\ga}{\ka_2}P-\dfrac{\ka_1+(\ga/2+2)\ka_2}{\ka_2}b.
\end{aligned}
\]
We use the above equation to simplify~\eqref{eq:couple2} into
\[
s_m(\del)c+c_m(\del)d
=\dfrac{3\ga}{4\ka_2\sinh\del}P
+\dfrac{(\ga+2)\ov{\ka}-4\ka_2}{4\ka_2\sinh\del}b.
\]
This gives~\eqref{eq:qccouple}$_1$, which together with~\eqref{eq:couple1} yields
~\eqref{eq:qccouple}$_2$.
\section{Derivation of Equation~\eqref{eq:qqc}}
To derive~\eqref{eq:qqc}, we firstly substitute the expression of $u_j$ into~\eqref{eq: qqc'} and obtain
\begin{align*}
&\quad\bigl((\ka_1+2\ka_2)(c_m(\del)-c_{m-1}(\del))
+\ka_2(c_{m+1}(\del)-c_{m-1}(\del))\bigr)c\\
&\quad+\bigl((\ka_1+2\ka_2)(s_m(\del)-s_{m-1}(\del))
+\ka_2(s_{m+1}(\del)-s_{m-1}(\del))\bigr)d=-P-\ov{\ka}b,
\end{align*}
and
\begin{align*}
&\quad\Bigl\{-(\ka_1+2\ka_2)(c_m(\del)-c_{m-1}(\del))\\
&\qquad+\ka_1(c_{m+1}(\del)-c_m(\del))+\ka_2(c_{m+2}(\del)-c_m(\del))\Bigr\}c\\
&\quad+\Bigl\{-(\ka_1+2\ka_2)(s_m(\del)-s_{m-1}(\del))\\
&\qquad\quad+\ka_1(s_{m+1}(\del)-s_m(\del))
+\ka_2(s_{m+2}(\del)-s_m(\del))\Bigr\}d=0.
\end{align*}
Using~\eqref{eq:pararela1}, we obtain
\begin{align}
(\ka_1+2\ka_2)(c_m(\del)&-c_{m-1}(\del))+\ka_2(c_{m+1}(\del)-c_{m-1}(\del))\nn\\
&=\Lr{\ka_1+2\ka_2+2(\cosh\del+1)\ka_2}(\cosh\del-1)c_{m-1}(\del)\nn\\
&\quad+(\ka_1+2\ka_2+2\ka_2\cosh\del)s_1(\del)s_{m-1}(\del)\nn\\
&=-\ov{\ka}c_{m-1}(\del).\label{eq:rela3}
\end{align}
Proceeding along the same line that leads to the above identity, we have
\[
(\ka_1+2\ka_2)(s_m(\del)-s_{m-1}(\del))
+\ka_2(s_{m+1}(\del)-s_{m-1}(\del))=-\ov{\ka}s_{m-1}(\del).
\]
Using the above two equations, we reshape the first equation of~\eqref{eq: qqc'} into~\eqref{eq:qqc}$_1$.

Using~\eqref{eq:rela3}, we write
\begin{align*}
&\quad-(\ka_1+2\ka_2)(c_m(\del)-c_{m-1}(\del))+\ka_1(c_{m+1}(\del)-c_m(\del))
+\ka_2(c_{m+2}(\del)-c_m(\del))\\
&=-\bigl\{(\ka_1+2\ka_2)(c_m(\del)-c_{m-1}(\del))
+\ka_2(c_{m+1}(\del)-c_m(\del))\bigr\}\\
&\quad+\bigl\{
(\ka_1+\ka_2)(c_{m+1}(\del)-c_m(\del))
+\ka_2(c_{m+2}(\del)-c_{m-1}(\del))\bigr\}\\
&=\ov{\ka}c_{m-1}(\del),
\end{align*}
where we have used~\eqref{eq:rela1} with $k=m+1$ in the last step.

Proceeding along the same line that leads to the above identity, we obtain
\[
-(\ka_1+2\ka_2)(s_m(\del)-s_{m-1}(\del))+\ka_1(s_{m+1}(\del)-s_m(\del))
+\ka_2(s_{m+2}(\del)-s_m(\del))=\ov{\ka}s_{m-1}(\del).
\]
Combining the above two equations, we obtain~\eqref{eq:qqc}$_2$.
\section{Derivation of~\eqref{eq:fqccouple}}
To derive~\eqref{eq:fqccouple}, we substitute the expression of $u_j$ into~\eqref{eq: fb'} and obtain
\[
\left\{
\begin{aligned}
\bigl(c_{m+1}(\del)-c_m(\del)\bigr)c&+\bigl(s_{m+1}(\del)-s_m(\del)\bigr)d
=-P/\ov{\ka},\\
\bigl(\ka_1(c_{m+2}(\del)-c_{m+1}(\del))
&+\ka_2(c_{m+3}(\del)-c_{m+1}(\del))\bigr)c\\
+\bigl(\ka_1(s_{m+2}(\del)-s_{m+1}(\del))
&+\ka_2(s_{m+3}(\del)-s_{m+1}(\del))\bigr)d\\
=-(\ka_1+2\ka_2)(P/\ov{\ka}+b).
\end{aligned}
\right.
\]
Proceeding along the same line that leads to~\eqref{eq:rela3}, we obtain
\begin{align*}
\ka_1\Lr{c_{m+2}(\del)-c_{m+1}(\del)}+\ka_2\Lr{c_{m+3}(\del)-c_{m+1}(\del)}
&=-2\ka_2s_1(\del)s_{m+1}(\del),\\
\ka_1\Lr{s_{m+2}(\del)-s_{m+1}(\del)}+\ka_2\Lr{s_{m+3}(\del)-s_{m+1}(\del)}
&=-2\ka_2s_1(\del)c_{m+1}(\del).
\end{align*}
We write the second equation of the coupling conditions as
\[
s_{m+1}(\del)c+c_{m+1}(\del)d
=\dfrac{\ka_1+2\ka_2}{2\ka_2\sinh\del}(P/\bar{\ka}+b)=-\coth\del(P/\bar{\ka}+b).
\]
This gives~\eqref{eq:fqccouple}$_1$.

In addition, we can write the first equation of the coupling equations as
\[
\bigl(c_{m+1}(\del)c+s_{m+1}(\del)d\bigr)(1-\cosh\del)
+\bigl(s_{m+1}(\del)c+c_{m+1}(\del)d\bigr)\sinh\del=-P/\bar{\ka}-b,
\]
which together with~\eqref{eq:fqccouple}$_1$ implies~\eqref{eq:fqccouple}$_2$.
\end{appendix}
\bibliographystyle{siam}
\bibliography{qcbif0}
\end{document}